\documentclass[twoside,11pt]{amsart}
\usepackage{amsmath,latexsym,amssymb,times,mathptm,amsfonts,enumerate,amsthm}
\usepackage{amsthm}
\usepackage{tikz}

\usepackage{pst-plot}
\usepackage{pst-node,pst-text,pst-3d,pstricks}

\def\proof{\noindent{\bf{Proof.} }}
\def\sqr#1#2{{\vcenter{\hrule height.#2pt
        \hbox{\vrule width.#2pt height#1pt \kern#1pt
                \vrule width.#2pt}
        \hrule height.#2pt}}}

\newtheorem{theorem}{Theorem}[section]

\newtheorem{lemma}[theorem]{Lemma}
\newtheorem{proposition}[theorem]{Proposition}

\newtheorem{definition}[theorem]{Definition}
\newtheorem{remark}[theorem]{Remark}

\newtheorem{example}[theorem]{Example}

\newcommand{\m}{\mathfrak{m}}

\newcommand{\ol}[1]{\ensuremath{\overline{#1}}}

\newcommand{\la}{\longrightarrow}
\newtheorem{construction}[theorem]{Construction}
\numberwithin{equation}{section}

\newcommand{\ds}{\displaystyle}

\newcommand{\rees}[1]{\ensuremath{\mathcal{R}{(#1)}}}
\newcommand{\assg}[1]{\ensuremath{{\rm gr}_{R}{(#1)}}}

\newcommand{\til}[1]{\ensuremath{\widetilde{#1}}}

\begin{document}
\baselineskip=16pt

\title{Rees algebras of square-free monomial ideals}
\author[L. Fouli]{Louiza Fouli}
\address{Department of Mathematical Sciences,  New Mexico State University,  Las Cruces, New Mexico 88003, USA}
\email{lfouli@math.nmsu.edu}

\author[K. N. Lin]{Kuei-Nuan Lin}
\address{Department of Mathematics, University of California, Riverside, 900 University Ave.,
Riverside, CA 92521, USA}
\email{linkuei@ucr.edu}


\begin{abstract} 
We study the defining equations of the Rees algebras of  square-free monomial ideals in a polynomial ring over a field. We propose the construction of a graph, namely the generator graph of a monomial ideal, where the monomial generators serve as vertices for the graph.
When $I$ is a square-free monomial ideal  such that each connected component of the generator graph of $I$  has at most $5$ vertices then $I$ has relation type at most $3$.
In general, we establish the defining equations of the Rees algebra in this case and give a combinatorial interpretation of them. Furthermore, we provide new classes of ideals of linear type.  We show that when $I$ is a square-free monomial ideal  and the generator graph of $I$ is the graph of a disjoint union of trees and graphs with a unique odd cycle, then $I$ is an ideal of linear type. 
\end{abstract}

\maketitle

\section{Introduction}

The main  problem of interest in this article is the question of determining  the defining equations
of the Rees algebra of a square-free monomial ideal in a polynomial ring over a field.
Let $R$ be a Noetherian ring and let $I$ be an ideal. The {\it Rees algebra}, $\rees{I}$, is  defined to be the graded algebra $\rees{I}=R[It]=\underset{i \geq 0} \oplus I^{i}t^{i} \subset R[t]$, where $t$ is an indeterminate. 
The Rees algebra of an ideal encodes many of the analytic properties of the variety defined by $I$ and it is the algebraic realization of the blowup of a variety along a subvariety.  The blowup  of ${\rm Spec}(R)$ along $V(I)$  is the projective spectrum of the Rees algebra, $\rees{I}$, of $I$. This  important construction is the main tool  in the resolution of singularities of an algebraic variety.

From the algebraic point of view the Rees algebra of an ideal facilitates the study of the asymptotic behavior of the ideal and it is essential in computing the integral closure of powers of ideals. The Rees algebra of an ideal $I$ in a Noetherian ring can be realized as a quotient of a polynomial ring and hence once the defining ideal of  $\rees{I}$ is determined, it is straightforward to compute the integral closure, $\ol{\rees{I}}$. It is well known that   $\ol{\rees{I}}= \underset{i \geq 0} \oplus \ol{I^i}t^i$ and one obtains $\ol{I^{i}}=[\ol{\rees{I}}]_{i}$.

 We consider the following construction for the Rees algebra. 
Let $R$ be a Noetherian ring and let $I$ be an $R$-ideal. Let $f_1, \ldots, f_n$ be a minimal generating set for $I$. Consider the polynomial ring $S=R[T_1, \ldots, T_n]$, where $T_1, \ldots, T_n$ are indeterminates. Then there is a natural map $\phi: S \longrightarrow \rees{I}=R[It]$ that sends $T_i$ to $f_it$. Let $J=\ker \phi$ be the defining ideal of $\rees{I}$. Then $\rees{I} \simeq S/J$ and  $J= \bigoplus \limits_{i=1}^{\infty}J_i $ is a graded ideal. A minimal generating set for $J$ is often referred to as the {\it defining equations} of the Rees algebra.
Also, $J_1$ is known as the ideal of linear relations as it is the defining ideal of the {\it symmetric algebra}, ${\rm Sym}(I)$, of $I$ and it is generated by linear forms in the $T_i$. The generators of $J_1$ arise from the first syzygies of $I$. When $J=J_1$ then $I$ is called an ideal of {\it linear type}. In this case $\rees{I}\simeq{\rm Sym}(I)$ and the defining ideal of $\rees{I}$ is well understood.

The problem of finding an explicit description of the defining equations of $\rees{I}$ has been addressed  by many authors, see for example \cite{Buse, CHTV, GerGim, GerGimHar, Gim, GimLor, Ha1, Ha2, KPU, MorTh, MorUlr, SimTrVal, Vas}.  In general, finding the generators for the defining ideal of the Rees algebra of an ideal is a difficult problem to address.

In Section 2, we develop a series of lemmas that allow us to determine conditions under which an element in the defining ideal of the Rees algebra becomes redundant and hence not needed in the defining equations. One of our essential tools is the result of Taylor, where she determines a non-minimal generating set for the defining ideal of the Rees algebra of a monomial ideal \cite{Tay}. One of our main results in this section give an explicit description of the defining equations of the Rees algebra of an ideal generated by $n \leq 5$ square-free monomials, Theorems~\ref{u gens}, ~\ref{5 gens linear}. 
Recall that the {\it relation type} of an ideal $I$ is defined to be ${\rm rt}(I)=\min \{s \mid J=\bigoplus \limits_{i=1}^{s}J_i\}$, where $J$ is the defining ideal of $\rees{I}$. In other words, the relation type of $I$ is the largest degree (in the $T_i$) of any minimal generator of $J$.  
In particular, when $I$ is of linear type it is of relation type
1.  The relation type of an ideal has been explored in various articles,
see for instance  \cite{AGH, Huc, Lai, Tru, Wa}. 
In Theorem~\ref{u gens} we establish that for a square-free monomial ideal generated by $n \leq 5$ monomials, the relation type is at most $n-2$. We also provide a class of examples of square-free monomial ideals generated by $n>4$ monomials, Example~\ref{rel type ex}, for which the relation type is at least $2n-7$. Notice that when $n=5$ then $2n-7=n-2$, and hence showing that the bound in Theorem~\ref{u gens} is attained. We remark that we do not impose any restrictions on the degrees of the monomial generators of the square-free monomial ideals under consideration.

We construct the {\it generator graph}  associated to monomial ideals, which is also known as the line graph of the hypergraph. The monomial generators are the vertices for this graph and the edges are determined by the greatest common divisor among the generators, see Construction~\ref{graph til}.  This construction enables us to extend Theorem~\ref{u gens} and allow more freedom on the number of generators of the ideal. We show that when $I$ is a square-free monomial ideal such that its generator graph is the disjoint union of graphs with at most $5$ vertices, then ${\rm rt}(I) \leq 3$, Theorem~\ref{what survives}. Furthermore, we show that the non-linear part of the defining equations of the Rees algebra, $\rees{I}$, of $I$ arise from even closed walks of the generator graph, Theorem~\ref{what survives}.

Recall that for a Noetherian local ring $(R, \m)$ and $I$ an $R$-ideal,  the special fiber ring, $\mathcal{F}(I)$, of $I$ is defined to be $\mathcal{F}(I)=\assg{I} \otimes R/\m$, where $\assg{I}=R[It]/IR[It]= \bigoplus \limits_{i=0}^{\infty} I^i/I^{i+1}$. When $R$ is a polynomial ring over a field $k$ and $I$ is an $R$-ideal then $\mathcal{F}(I)$ is defined to be $k[f_1, \ldots, f_n]$, where $f_1,\ldots, f_n$ is a minimal generating set for $I$. Therefore, there is a homomorphism  $\psi: k[T_1, \ldots, T_n] \rightarrow \mathcal{F}(I)$ that sends $T_i$ to $f_i$. Let $H=\ker \psi$ and $\mathcal{F}(I)\simeq k[T_1, \ldots, T_n]/H$. Let $S=R[T_1, \ldots, T_n]$ and let $J$ be the defining ideal of the Rees algebra, $\rees{I}\simeq S/J$, of $I$. An ideal $I$ is called an ideal of {\it fiber type} if the defining ideal $J$ of $\rees{I}$ is obtained by the linear relations and the defining equations for the special fiber ring, i.e. $J=SJ_1+SH$. In the case of edge ideals, Villarreal showed that they are all ideals of fiber type, \cite[Theorem~3.1]{Vill}.  Moreover, Villarreal gave an explicit description of the defining equations of the Rees algebra of any square-free monomial ideal generated in degree $2$, \cite[Theorem~3.1]{Vill}.
It is also worth noting that Villarreal exhibited an example to show that his techniques do not extend for monomial ideals generated in degree higher than $2$, \cite[Example~3.1]{Vill}. His example is a square-free monomial ideal generated in degree $3$ that is not of fiber type. We remark that in general, a square-free monomial ideal generated in degree greater or equal to$3$ is not of fiber type. In Example~\ref{Vil ex} we give a concrete explanation of how one obtains the defining equations for the ideal in \cite[Example~3.1]{Vill}.

In Section~\ref{lin type}, we concentrate on the special class of ideals of linear type. The first well known class of ideals of linear type are complete intersection ideals, \cite{Micali}. Huneke and Valla showed that ideals generated by $d$-sequences are also ideals of linear type, \cite{Hu1}, \cite{Valla}. Later on, Herzog-Simis-Vasconcelos and Herzog-Vasconcelos-Villarreal used strongly Cohen-Macaulay ideals and certain conditions on local numbers of generators to describe other classes of ideals of linear type, \cite{HSV1, HSV2, HVV}. Huneke also showed that the ideal of all $(n-1)\times (n-1)$ minors of a generic matrix in $R=\mathbb{Z}[X_{i,j}]$ is also an ideal of linear type, \cite{Hu2}. In the case of square--free monomial ideals generated in degree $2$, Villarreal gave a complete characterization of ideals of linear type. More precisely, he showed that such an ideal $I$ is of linear type if and only if the graph of $I$  is the disjoint union of graphs of trees or graphs that have a unique odd cycle, \cite[Corollary~3.2]{Vill}.
In this work we find new classes of square-free monomial
ideals that are of linear type. We prove that when $I$ is an ideal generated by square-free monomials and the generator graph of $I$ is a
disjoint union of graphs of trees and graphs with a unique odd cycle then $I$ is an ideal of linear type, Theorem~\ref{union graph linear type}. Our results have a similar flavor as those of Villarreal in \cite{Vill}. In some sense these results extend the work of Villarreal, even though we do not fully recover his results. Nonetheless, our techniques allow us to consider square-free monomial ideals without any restrictions on the degrees of the generators.

\section{The defining equations of the Rees algebra}

Let $R$ be a polynomial ring over a field and let $I$ be a  monomial ideal in $R$. Let $f_1, \ldots, f_n$ be a minimal monomial generating set of $I$ and
let $\rees{I}$ denote the Rees algebra of $I$. Then $\rees{I}=R[f_1t, \ldots, f_nt] \subset R[t]$ and there is an epimorphism $$\phi: S=R[T_1, \ldots, T_n] \la \rees{I}$$
induced by $\phi(T_i)=f_it$. Let $J=\ker \phi$. We note that $J=\bigoplus \limits_{i=1}^{\infty}J_i$ is a graded ideal of $S$. As mentioned in the Introduction, a minimal generating set for the ideal $J$ is referred to as the defining equations of the Rees algebra of $I$.

\begin{definition}\label{sequences}{\rm
Let $I$ be a monomial ideal in a polynomial ring $R$ over a field $k$. Let $f_{1}, \ldots, f_{n}$ be a minimal monomial generating set of $I$.
 Let $\mathcal{I}_s$ denote the set of all non-decreasing sequences of integers $\alpha=(i_1, \ldots, i_s)\subset\{1,\ldots.,n\}$ of length $s$. Then $f_{\alpha}=f_{i_1}\cdots f_{i_s}$ is the corresponding product of monomials in $I$.  We also let $T_{\alpha}=T_{i_1} \cdots T_{i_s}$  be the corresponding product of indeterminates in $S=R[T_1, \ldots, T_n]$. For every $\alpha, \beta \in \mathcal{I}_s$ we consider the binomial $T_{\alpha, \beta}=\ds{\frac{f_{\beta}}{\gcd(f_{\alpha}, f_{\beta})}T_{\alpha}-\frac{f_{\alpha}}{\gcd(f_{\alpha},f_{\beta})}T_{\beta}}$.}
\end{definition}

 The following is in \cite{Tay} and we cite it here for ease of reference.

\begin{theorem}$($\cite{Tay}$)$ \label{Taylor}
Let $R$ be a polynomial ring over a field and let $I$ be a monomial ideal in $R$.  Let $f_{1}, \ldots, f_{n}$ be a minimal monomial generating set of $I$ and let $\rees{I}$ be the Rees algebra of $I$. 
Then $\rees{I}\simeq S/J$, where $S=R[T_1, \ldots, T_n]$ is a polynomial ring, $T_1, \ldots, T_n$ are indeterminates, and $J=SJ_1 + S \cdot (\bigcup \limits_{i=2}^{\infty}J_i)$ such that $J_s=\{ T_{\alpha, \beta} \mid \mbox{ for }\alpha, \beta \in \mathcal{I}_{s}\}$. 
\end{theorem}

The following remark states some properties for the greatest common divisor among square-free monomials. The proofs are omitted as they are elementary.

\begin{remark}\label{more gcd prop}
Let $a,b,c,d, b_1, \ldots, b_s$ be square-free monomials  in a polynomial ring $R$ over a field. Let $n,m,l,r $ be positive integers. Then
\begin{enumerate}[$($a$)$]
\item $\gcd (a^n,b^m)=\gcd(a,b)^{\min\{m,n\}}$,
\item $\gcd(a^{n}, b_1 \cdots b_s)=\gcd(a^{s}, b_1 \cdots b_s)$ for all $ n \geq s$,
\item $\gcd(a^n,b^mc^l)=\ds{\frac{\gcd(a,b)^{\min \{n,m\}}\gcd(a,c)^{\min\{n,l\}}}{C}}$, for some $C\in R$.

\end{enumerate}
\end{remark}

The main goal in this section is to determine conditions under which for sequences $\alpha, \beta$ of length $s \geq 2$ we have $T_{\alpha, \beta} \in SJ_1+S \cdot (\bigcup \limits_{i=2}^{s-1}J_i)$. In other words, we are interested in finding conditions on the sequences $\alpha, \beta$ such that the generator $T_{\alpha, \beta}$ is  redundant in the defining ideal of the Rees algebra. 
The following two lemmas follow directly from Definition~\ref{sequences}.

\begin{lemma}\label{L3}
Let $R$ be a polynomial ring over a field and let $I$ be a  monomial ideal in $R$.  Let $\alpha=(a_1, \ldots, a_s)$, $\beta=(b_1, \ldots, b_s) \in \mathcal{I}_s$  be two sequences of length $s \geq 2$, where $\mathcal{I}_s$ is as in Definition~\ref{sequences}.  Suppose that $a_i=b_j$ for some $i$ and some $j$. Then $T_{\alpha, \beta}  \in SJ_1 + S \cdot (\bigcup \limits_{i=2}^{s-1}J_i)$.\end{lemma}

\proof
Without loss of generality we may assume that $a_1=b_1$. Let $\alpha_1=(a_2, \ldots, a_s)$ and $\beta_1=(b_2, \ldots, b_s)$. Then $\gcd(f_{\alpha},f_{\beta})=f_{a_1}\gcd(f_{\alpha_1}, f_{\beta_1})=f_{b_1}\gcd({f_{\alpha_1}},{f_{\beta_1}})$. Notice that $f_{a_1}=f_{b_1}$ and $T_{a_1}=T_{b_1}$.
Then 
\begin{eqnarray*}
T_{\alpha, \beta}&=&\ds{\frac{ f_{\beta}}{\gcd(f_{\alpha}, f_{\beta})}T_{\alpha}-\frac{f_{\alpha}}{\gcd(f_{\alpha}, f_{\beta})}T_{\beta}}
= T_{a_1}[\ds{\frac{{f_{\beta_1}}}{\gcd({f_{\alpha_1}}, {f_{\beta_1}})}{T_{\alpha_1}}-\frac{{f_{\alpha_1}}}{\gcd({f_{\alpha_1}},{f_{\beta_1}})}{T_{\beta_1}}}]\\
&=&T_{a_1}{[T_{\alpha_1, \beta_1}]}\in SJ_1 + S \cdot (\bigcup \limits_{i=2}^{s-1}J_i). \hspace{8.5 cm} \qed
\end{eqnarray*}

\begin{lemma}\label{L4}
Let $R$ be a polynomial ring over a field and let $I$ be a monomial ideal in $R$.  Let $\alpha$, $\beta \in \mathcal{I}_s$  be two sequences of length $s \geq 2$, where $\mathcal{I}_s$ is as in Definition~\ref{sequences}. Suppose that $\alpha=(\alpha_1, \ldots, \alpha_1)$ and $\beta=(\beta_1, \ldots, \beta_1)$, where $\alpha_1, \beta_1 \in \mathcal{I}_{m}$ and $m \leq s$. Then $T_{\alpha, \beta}  \in SJ_1 + S \cdot (\bigcup \limits_{i=2}^{m}J_i)$.
\end{lemma}

\proof 
First we observe that $s$ is a multiple of $m$. Let $l$ be an integer such that $s=lm$.
Notice that $f_{\alpha}=f_{\alpha_1}^{l}$ and $f_{\beta}=f_{\beta_1}^{l}$.  Also, by Remark~\ref{more gcd prop}~(a) we have that $\gcd(f_{\alpha}, f_{\beta})=\gcd(f_{\alpha_1},f_{\beta_1})^l$.  Let $a=\ds{\frac{f_{\beta_1} T_{\alpha_1}}{\gcd(f_{\alpha_1},f_{\beta_1})}}$, $b=\ds{\frac{f_{\alpha_1} T_{\beta_1}}{\gcd(f_{\alpha_1},f_{\beta_1})}}$, and let $A=(a^{l-1}+a^{l-2}b+\ldots +ab^{l-2}+b^{l-1})$. 
Notice that $a,b \in S$ and thus $A \in S$. Then 
\begin{eqnarray*}
T_{\alpha, \beta}&=& \ds{\frac{ f_{\beta}}{\gcd(f_{\alpha}, f_{\beta})}T_{\alpha}-\frac{f_{\alpha}}{\gcd(f_{\alpha}, f_{\beta})}T_{\beta}}=a^l-b^l\\
&=&(a-b)(a^{l-1}+a^{l-2}b+\ldots +ab^{l-2}+b^{l-1})=T_{a_1, b_1}A  \in SJ_1+S \cdot (\bigcup \limits_{i=2}^{m}J_{i}). \hfill \qed
\end{eqnarray*}

We observe the following properties for greatest common divisors among monomials. Again, we omit the proofs as they are elementary.
\begin{remark} \label{gcd properties}
Let $a, b, c, d,e,f$ be monomials in a polynomial ring $R$ over a field. 
\begin{enumerate}[$($a$)$]
\item ${\gcd}(a, bc)=\ds{\frac{{\gcd}(a,b){\gcd}(a,c)}{C}}$, for some $C \in R$.
\item Suppose that $\gcd(a,c)=1$. Then $\gcd(ab,c)=\gcd(b,c)$.
\end{enumerate}
\end{remark}

The following Lemma plays an important role in the rest of this article as it allows us to show that certain expressions $T_{\alpha, \beta}$ are redundant in the defining ideal of the Rees algebra.

\begin{lemma}\label{Tab split extended}
Let $R$ be a polynomial ring over a field and let $I$ be a  monomial ideal in $R$.  Let $\alpha$, $\beta \in \mathcal{I}_s$  be two sequences of length $s \geq 2$, where $\mathcal{I}_s$ is as in Definition~\ref{sequences}. Suppose
$\alpha=(\alpha_{1},\ldots,\alpha_{m})$ and $\beta=(\beta_{1},\ldots,\beta_{m})$
with $\alpha_{i}$, $\beta_{i}\in I_{s_{i}}$ and $s_{1}+\ldots+s_{m}=s$.
Suppose that  for all $1 \leq i \leq m$, there exist $C_{i}\in R$ such that \[
\gcd(f_{\alpha},f_{\beta})=\ds{\frac{\gcd(\prod \limits_{j=1}^{i-1}f_{\alpha_{j}},f_{\beta})\gcd(\prod \limits_{k=i}^{m}f_{\alpha_{k}},\prod \limits_{k=i+1}^{m}f_{\beta_{k}})\gcd(f_{\alpha_{i}},f_{\beta_{i}})}{C_{i}}},\]
where empty products are taken to be $1$.
Then 
$$T_{\alpha,\beta}=\ds{\sum \limits_{i=1}^{m}\left[A_{i}\prod_{j=i+1}^{m}T_{\alpha_{j}}\prod_{k=1}^{i-1}T_{\beta_{k}}\right]T_{\alpha_{i},\beta_{i}} }\in SJ_1 + S \cdot (\bigcup \limits_{i=2}^{r}J_i),$$
where $r=\max\{s_1, \ldots, s_m\}$ and $A_i \in R$ for all $i$.
\end{lemma}

\proof
Notice that $f_{\alpha}=\prod \limits_{i=1}^{m}f_{\alpha_{i}}$ and $f_{\beta}=\prod \limits_{i=1}^{m}f_{\beta_{i}}$. Then we have \begin{eqnarray*}
T_{\alpha,\beta} & = & \frac{f_{\beta}}{\gcd(f_{\alpha},f_{\beta})}T_{\alpha}-\frac{f_{\alpha}}{\gcd(f_{\alpha},f_{\beta})}T_{\beta}\\
 & = &\ds{ \sum_{i=1}^{m}\left[\frac{C_{i}\prod \limits_{j=1}^{i-1}f_{\alpha_{j}}\prod \limits_{k=i+1}^{m}f_{\beta_{k}}\prod \limits_{j=i+1}^{m}T_{\alpha_{j}}\prod \limits_{k=1}^{i-1}T_{\beta_{k}}}{\gcd(\prod \limits_{j=1}^{i-1}f_{\alpha_{j}},f_{\beta})\gcd(\prod \limits_{k=i}^{m}f_{\alpha_{k}},\prod \limits_{k=i+1}^{m}f_{\beta_{k}})}  \left(\frac{f_{\beta_{i}}}{\gcd(f_{\alpha_{i}},f_{\beta_{i}})}T_{\alpha_{i}}-\frac{f_{\alpha_{i}}}{\gcd(f_{\alpha_{i}},f_{\beta_{i}})}T_{\beta_{i}}\right)\right]}. \end{eqnarray*}
Finally, we note that $A_i=\ds{\frac{C_{i}\prod \limits_{j=1}^{i-1}f_{\alpha_{j}}\prod \limits_{k=i+1}^{m}f_{\beta_{k}}}{\gcd(\prod \limits_{j=1}^{i-1}f_{\alpha_{j}},f_{\beta})\gcd(\prod \limits_{k=i}^{m}f_{\alpha_{k}},\prod \limits_{k=i+1}^{m}f_{\beta_{k}})}} \in R$. \qed

\bigskip

Next we give a list of conditions that when satisfied by two sequences $\alpha, \beta$ then the generator $T_{\alpha, \beta}$ is not a minimal generator in the defining ideal of the Rees algebra of the corresponding ideal.

\begin{proposition}\label{L2}
Let $R$ be a polynomial ring over a field and let $I$ be a monomial ideal in $R$.  Let $\alpha$, $\beta \in \mathcal{I}_s$  be two sequences of length $s \geq 2$, where $\mathcal{I}_s$ is as in Definition~\ref{sequences}.  Let  $\alpha=(a_1, \ldots, a_s), \beta=(b_1, \ldots, b_s)$. Suppose that after some reordering there exist integers $k, l$ with $1 \leq k,l \leq s-1$ such that $\gcd (f_{a_{i}}, f_{b_j})=1$ for every $1 \leq i \leq l$ and every $k+1 \leq j \leq s$. We further assume that $\gcd (f_{a_u},f_{b_{v}})=1$, for every $l+1 \leq u \leq s$ and every $1 \leq v \leq k$. Then $T_{\alpha, \beta} \in SJ_1 + S \cdot (\bigcup \limits_{i=2}^{r}J_i)$, where $r=\max \{k, | k-l|, s-k, s-l\}$. 
\end{proposition}

\proof
Suppose that $l=k$ and  write $\alpha=(\alpha_{1},\alpha_{2})$ and $\beta=(\beta_{1},\beta_{2})$,
where $\alpha_{1}=(a_{1},\ldots,a_{k})$, $\alpha_{2}=(a_{k+1},\ldots,a_{s})$,
$\beta_{1}=(b_{1},\ldots,b_{k})$, and  $\beta_{2}=(b_{k+1},\ldots,b_{s})$.
By our assumptions, we have $\gcd(f_{\alpha_{i}},f_{\beta_{j}})=1$
for $i\neq j$. Hence by Remark~\ref{gcd properties}, we have  \begin{eqnarray*}
\gcd(f_{\alpha},f_{\beta}) & = & \frac{\gcd(f_{\alpha},f_{\beta_{1}})\gcd(f_{\alpha},f_{\beta_{2}})}{C_{1}}
 =  \frac{\gcd(f_{\alpha_{1}},f_{\beta_{1}})\gcd(f_{\alpha},f_{\beta_{2}})}{C_{1}},\\
\gcd(f_{\alpha},f_{\beta}) & = & \frac{\gcd(f_{\alpha_{1}},f_{\beta})\gcd(f_{\alpha_{2}},f_{\beta})}{C_{2}}=  \frac{\gcd(f_{\alpha_{1}},f_{\beta})\gcd(f_{\alpha_{2}},f_{\beta_{2}})}{C_{2}}.\end{eqnarray*}
The result follows from Lemma~\ref{Tab split extended} with $m=2$.

Without loss of generality we may now assume that $l<k$. We write
$\alpha=(\alpha_{1},\alpha_{2},\alpha_{3})$ and $\beta=(\beta_{1},\beta_{2},\beta_{3})$,
with $\alpha_{1}=(a_{1},\ldots,a_{l})$, $\alpha_{2}=(a_{l+1},\ldots,a_{k})$,
$\alpha_{3}=(a_{k+1},\ldots,a_{s})$, $\beta_{1}=(b_{1},\ldots,b_{l})$,
$\beta_{2}=(b_{l+1},\ldots,b_{k})$, and $\beta_{3}=(b_{k+1},\ldots,b_{s})$.
Then $\gcd(f_{\alpha_{1}},f_{\beta_{3}})=1$
and $\gcd(f_{\alpha_{i}},f_{\beta_{j}})=1$ for $i=2,3$ and $j=1,2$.
Hence by Remark~\ref{gcd properties} we have

\begin{eqnarray*}
\gcd(f_{\alpha},f_{\beta}) & = & \frac{\gcd(f_{\alpha},f_{\beta_{1}})\gcd(f_{\alpha},f_{\beta_{2}}f_{\beta_{3}})}{C_{1}} =  \frac{\gcd(f_{\alpha_{1}},f_{\beta_{1}})\gcd(f_{\alpha},f_{\beta_{2}}f_{\beta_{3}})}{C_{1}},\\
 \gcd(f_{\alpha},f_{\beta}) & = & \frac{\gcd(f_{\alpha_{1}},f_{\beta})\gcd(f_{\alpha_{2}}f_{\alpha_{3}},f_{\beta})}{C_{2}}\\
 & = & \frac{\gcd(f_{\alpha_{1}},f_{\beta})\gcd(f_{\alpha_{2}}f_{\alpha_{3}},f_{\beta_{3}})\gcd(f_{\alpha_{2}}f_{\alpha_{3}},f_{\beta_{1}}f_{\beta_{2}})}{C_{2}'}\\
 & = & \frac{\gcd(f_{\alpha_{1}},f_{\beta})\gcd(f_{\alpha_{2}}f_{\alpha_{3}},f_{\beta_{3}})\gcd(f_{\alpha_{2}}f_{\alpha_{3}},f_{\beta_{2}})}{C_{2}'}\\
 & = & \frac{\gcd(f_{\alpha_{1}},f_{\beta})\gcd(f_{\alpha_{2}}f_{\alpha_{3}},f_{\beta_{3}})\gcd(f_{\alpha_{2}},f_{\beta_{2}})}{C_{2}'}, \mbox{ and}\\
\gcd(f_{\alpha},f_{\beta}) & = & \frac{\gcd(f_{\alpha_{1}}f_{\alpha_{2}},f_{\beta})\gcd(f_{\alpha_{3}},f_{\beta})}{C_{3}}=  \frac{\gcd(f_{\alpha_{1}}f_{\alpha_{2}},f_{\beta})\gcd(f_{\alpha_{3}},f_{\beta_{3}})}{C_{3}}.\end{eqnarray*}
Therefore we may apply Lemma~\ref{Tab split extended} with $m=3$.  \qed

The next three Lemmas deal with possible repetitions in the sequences $\alpha, \beta \in \mathcal{I}_s$ and how that affects the term $T_{\alpha, \beta}$.

\begin{lemma}\label{L5}
Let $R$ be a polynomial ring over a field and let $I$ be a square-free monomial ideal in $R$.  Let $\alpha$, $\beta \in \mathcal{I}_s$  be two sequences of length $s \geq 2$, where $\mathcal{I}_s$ is as in Definition~\ref{sequences}. Suppose that $\alpha=(a_1, \ldots, a_1)$ and  $\beta=(b_1, \ldots,  b_s)$. Then $T_{\alpha, \beta} \in SJ_1 + S \cdot (\bigcup \limits_{i=2}^{s-1}J_i)$.
\end{lemma}

\proof Notice that $f_{\alpha}=f_{a_1}^s$ and $f_{\beta}=f_{b_1} \cdots f_{b_s}$.  Let ${\beta_1}=(b_2, \ldots, b_s)$. We also note that $\gcd(f_{a_1}^s, f_{b_1})=\gcd(f_{a_1}, f_{b_1})$, by Remark~\ref{more gcd prop}~(a). Then 
\begin{eqnarray*}\gcd(f_{\alpha}, f_{\beta})&=&\ds{\frac{\gcd(f_{a_1}^{s}, f_{b_1})\gcd(f_{a_1}^s,{f_{\beta_1}})}{C}}=\ds{\frac{\gcd(f_{a_1}, f_{b_1})\gcd(f_{a_1}^{s-1},{f_{\beta_1}})}{C}},
\end{eqnarray*} for some $C \in R$, by Remark~\ref{more gcd prop}~(b) and (c) . 
Hence $T_{\alpha, \beta} \in SJ_1 + S \cdot (\bigcup \limits_{i=2}^{s-1}J_i)$, by Lemma~\ref{Tab split extended}. \qed

\begin{lemma}\label {L6}
Let $R$ be a polynomial ring over a field and let $I$ be a square-free monomial ideal in $R$. Let $\alpha$, $\beta \in \mathcal{I}_s$  be two sequences of length $s \geq 2$, where $\mathcal{I}_s$ is as in Definition~\ref{sequences}. Suppose that $\alpha=(a_1, \ldots, a_1, a_2, \ldots, a_2)$ and $\beta=(b_1, \ldots, b_1, b_2, \ldots, b_2)$. Then $T_{\alpha, \beta} \in SJ_1 + SJ_2$. \end{lemma}

\proof 
We will proceed by induction on $s$. If $s=2$ then there is nothing
to show. Suppose that $s>2$. Suppose that there are $l_{i}$ distinct copies
of $a_{i}$ in $\alpha$ and $k_{i}$ distinct copies of $b_{i}$ in $\beta$
for $i=1,2$. Then $f_{\alpha}=f_{a_{1}}^{l_{1}}f_{a_{2}}^{l_{2}}$
and $f_{\beta}=f_{b_{1}}^{k_{1}}f_{b_{2}}^{k_{2}}$. For simplicity
we write $\alpha=(a_{1}^{l_{1}},a_{2}^{l_{2}})$ and $\beta=(b_{1}^{k_{1}},b_{2}^{k_{2}})$.
If $k_{1}=k_{2}$ and $l_{1}=l_{2}$ then the result follows from
Lemma~\ref{L4}. Thus we may assume without loss of generality that $k_{1}>k_{2}$
and $l_{1}\geq l_{2}$.

We will show \begin{eqnarray*}
\gcd(f_{\alpha},f_{\beta}) & = & \frac{\gcd(f_{a_{1}}^{l_{1}-1}f_{a_{2}}^{l_{2}},f_{b_{1}}^{k_{1}}f_{b_{2}}^{k_{2}})\gcd(f_{a_{1}},f_{b_{1}})}{G}\\
 & = & \frac{\gcd(f_{a_{1}}^{l_{1}-1}f_{a_{2}}^{l_{2}},f_{b_{1}}^{k_{1}-1}f_{b_{2}}^{k_{2}})\gcd(f_{a_{1}}^{l_{1}}f_{a_{2}}^{l_{2}},f_{b_{1}})}{H}\end{eqnarray*}
 for some $G, H \in R$. Then $T_{\alpha, \beta}  \in SJ_1 \cup SJ_2$, by Lemma~\ref{Tab split extended}  and the induction. 

Notice that we have the following equalities: \begin{eqnarray*}
\gcd(f_{\alpha},f_{\beta}) & = & \frac{\gcd(f_{a_{1}}^{l_{1}-1}f_{a_{2}}^{l_{2}},f_{b_{1}}^{k_{1}}f_{b_{2}}^{k_{2}})\gcd(f_{a_{1}},f_{b_{1}}^{k_{1}}f_{b_{2}}^{k_{2}})}{C}\\
 & = & \frac{\gcd(f_{a_{1}}^{l_{1}-1}f_{a_{2}}^{l_{2}},f_{b_{1}}^{k_{1}}f_{b_{2}}^{k_{2}})\gcd(f_{a_{1}},f_{b_{1}}^{k_{1}-1}f_{b_{2}}^{k_{2}})\gcd(f_{a_{1}},f_{b_{1}})}{CD},\\
 \gcd(f_{\alpha},f_{\beta}) & = & \frac{\gcd(f_{a_{1}}^{l_{1}}f_{a_{2}}^{l_{2}},f_{b_{1}}^{k_{1}-1}f_{b_{2}}^{k_{2}})\gcd(f_{a_{1}}^{l_{1}}f_{a_{2}}^{l_{2}},f_{b_{1}})}{E}\\
 & = & \frac{\gcd(f_{a_{1}}^{l_{1}-1}f_{a_{2}}^{l_{2}},f_{b_{1}}^{k_{1}-1}f_{b_{2}}^{k_{2}})\gcd(f_{a_{1}},f_{b_{1}}^{k_{1}-1}f_{b_{2}}^{k_{2}})\gcd(f_{a_{1}}^{l_{1}}f_{a_{2}}^{l_{2}},f_{b_{1}})}{EF},\end{eqnarray*}
for some $C, D, E, F \in R$ by Remark~\ref{gcd properties}~(a). We will show that for any
integer $t$ and any variable $x \in R$ if $x^{t}\mid\gcd(f_{a_{1}},f_{b_{1}}^{k_{1}}f_{b_{2}}^{k_{2}})$,
then $x^{t}\mid CD$ and $x^{t}\mid EF$. This is equivalent
to showing that for any variable $x \in R$ and any integer $u$, if $x \mid\gcd(f_{a_{1}},f_{b_{1}}^{k_{1}}f_{b_{2}}^{k_{2}})$
and $x^{u}\mid\gcd(f_{\alpha},f_{\beta})$ then $x^{u}\mid\gcd(f_{a_{1}}^{l_{1}-1}f_{a_{2}}^{l_{2}},f_{b_{1}}^{k_{1}}f_{b_{2}}^{k_{2}})\gcd(f_{a_{1}},f_{b_{1}})$
and $x^{u}\mid\gcd(f_{a_{1}}^{l_{1}-1}f_{a_{2}}^{l_{2}},f_{b_{1}}^{k_{1}-1}f_{b_{2}}^{k_{2}})\gcd(f_{a_{1}}^{l_{1}}f_{a_{2}}^{l_{2}},f_{b_{1}})$.
This follows immediately since $s=l_{1}+l_{2}=k_{1}+k_{2}>l_{1}$ and $l_{1}\geq \frac{s}{2}>k_{2}$. \qed

\bigskip

\begin{lemma}\label {5gensub}
Let $R$ be a polynomial ring over a field and let $I$ be a square-free monomial ideal in $R$.  Let $\alpha$, $\beta \in \mathcal{I}_s$  be two sequences of length $s \geq 4$, where $\mathcal{I}_s$ is as in Definition~\ref{sequences}. Suppose that $\alpha=(a_{1},\ldots,a_{1},a_{2},\ldots,a_{2},a_{3},\ldots,a_{3})$ and $\beta=(b_{1},\ldots,b_{1},b_{2},\ldots,b_{2})$. Then $T_{\alpha,\beta}\in SJ_{1}+ S\cdot (J_{2}\cup J_{3})$. \end{lemma}

\proof 
We will show $T_{\alpha,\beta}\in SJ_1 + S \cdot (\bigcup \limits_{i=2}^{s-1}J_i)$ and  by
induction on $s \geq 4$, we conclude that $T_{\alpha,\beta}\in SJ_{1}+ S\cdot (J_{2}\cup J_{3})$. Suppose that there are  $l_{i}$ distinct copies of $a_{i}$ in $\alpha$ and $k_{i}$ distinct copies of $b_{i}$ in $\beta$. Then $l_{1}+l_{2}+l_{3}=k_{1}+k_{2}=s\geq4$.
Without lost of generality, we assume that $l_{1}\geq l_{2}\geq l_{3}\geq1$
and $k_{1}\geq k_{2}\geq1$. Since $s\geq4$, we have $l_{1}\geq2$
and $k_{1}\geq2$. We write $\alpha=(a_{1}^{l_{1}}, a_{2}^{l_{2}}, a_{3}^{l_{3}})$
and $\beta=(b_{1}^{k_{1}},b_{2}^{k_{2}})$. We then have three
possible scenarios. We claim the following:

\begin{enumerate}[$($i$)$]
\item If $l_{1}+l_{2}>k_{1}$, then  \begin{eqnarray*}
\gcd(f_{\alpha},f_{\beta}) & = & \gcd(f_{a_{1}}f_{a_{2}},f_{b_{1}}f_{b_{2}})\gcd(f_{\alpha_{1}},f_{\beta})/A\\
 & = & \gcd(f_{\alpha},f_{b_{1}}f_{b_{2}})\gcd(f_{\alpha_{1}},f_{\beta_{1}})/B,\end{eqnarray*}
 for some  $A, B \in R$, and $\alpha_{1}=(a_{1}^{l_{1}-1}, a_{2}^{l_{2}-1}, a_{3}^{l_{3}})$
and $\beta_{1}=(b_{1}^{k_{1}-1}, b_{2}^{k_{2}-1})$. 

\item If $l_{1}>k_{2}$, then  \begin{eqnarray*}
\gcd(f_{\alpha},f_{\beta}) & = & \gcd(f_{a_{1}},f_{b_{1}})\gcd(f_{\alpha_{1}},f_{\beta})/C\\
 & = & \gcd(f_{\alpha},f_{b_{1}})\gcd(f_{\alpha_{1}},f_{\beta_{1}})/D,\end{eqnarray*}
 for some  $C, D \in R$, and $\alpha_{1}=(a_{1}^{l_{1}-1}, a_{2}^{l_{2}}, a_{3}^{l_{3}})$
and $\beta_{1}=(b_{1}^{k_{1}-1}, b_{2}^{k_{2}})$. 

\item If $l_{1}+l_{2}\leq k_{1}$ and $l_{1}\leq k_{2}$, then \begin{eqnarray*}
\gcd(f_{\alpha},f_{\beta})&=&\gcd(f_{a_{1}}^{s/3}f_{a_{2}}^{s/3}f_{a_{3}}^{s/3},f_{b_{1}}^{2s/3}f_{b_{2}}^{s/3})\\
&=&(\gcd(f_{a_{1}}f_{a_{2}}f_{a_{3}},f_{b_{1}}^{2}f_{b_{2}}))^{s/3}. \end{eqnarray*}

\end{enumerate}
Once we obtain the above claims then the result will follow by Lemma \ref{Tab split extended}.

To establish claim (i) we notice that we have the following equalities:\begin{eqnarray*}
\gcd(f_{\alpha},f_{\beta}) & = & \gcd(f_{a_{1}}f_{a_{2}},f_{\beta})\gcd(f_{\alpha_{1}},f_{\beta})/E\\
 & = & \gcd(f_{a_{1}}f_{a_{2}},f_{b_{1}}f_{b_{2}})\gcd(f_{a_{1}}f_{a_{2}},f_{\beta_{1}})\gcd(f_{\alpha_{1}},f_{\beta})/EF,\\
 & = & \gcd(f_{\alpha},f_{b_{1}}f_{b_{2}})\gcd(f_{\alpha},f_{\beta_{1}})/G\\
 & = & \gcd(f_{\alpha},f_{b_{1}}f_{b_{2}})\gcd(f_{\alpha_{1}},f_{\beta_{1}})\gcd(f_{a_{1}}f_{a_{2}},f_{\beta_{1}})/GH,\end{eqnarray*}
for some $E, F, G, H \in R$. It is enough to
show that for any variable $x \in R$ and any positive integer $t$ if $x^{t}\mid\gcd(f_{a_{1}}f_{a_{2}},f_{\beta_{1}})$,
then $x^{t}\mid EF$ and $x^{t}\mid GH$. This is equivalent to showing that if $x^{u}\mid\gcd(f_{\alpha},f_{\beta})$
and $x\mid\gcd(f_{a_{1}}f_{a_{2}},f_{\beta_{1}})$, then $x^{u}\mid\gcd(f_{a_{1}}f_{a_{2}},f_{b_{1}}f_{b_{2}})\gcd(f_{\alpha_{1}},f_{\beta})$
and $x^{u}\mid\gcd(f_{\alpha},f_{b_{1}}f_{b_{2}})\gcd(f_{\alpha_{1}},f_{\beta_{1}})$
for any positive integer $u$. But this follows immediately since $l_{1}+l_{2}>k_{1}\geq k_{2}$
and $l_{1}+l_{2}<l_{1}+l_{2}+l_{3}=s=k_{1}+k_{2}$.

For claim (ii) we notice that we have the following equalities:\begin{eqnarray*}
\gcd(f_{\alpha},f_{\beta}) & = & \gcd(f_{a_{1}},f_{\beta})\gcd(f_{\alpha_{1}},f_{\beta})/K\\
 & = & \gcd(f_{a_{1}},f_{b_{1}})\gcd(f_{a_{1}},f_{\beta_{1}})\gcd(f_{\alpha_{1}},f_{\beta})/KL,\\
 & = & \gcd(f_{\alpha},f_{b_{1}})\gcd(f_{\alpha},f_{\beta_{1}})/M\\
 & = & \gcd(f_{\alpha},f_{b_{1}})\gcd(f_{\alpha_{1}},f_{\beta_{1}})\gcd(f_{a_{1}},f_{\beta_{1}})/MN,\end{eqnarray*}
for some $K, L, M, N \in R$. The  claim follows since  $l_{1}>k_{2}$.

Finally to establish the last claim we note that the assumption $l_{1}\leq k_{2}$ is equivalent to
$l_{2}+l_{3}\geq k_{1}$ and thus $l_{2}+l_{3}\geq k_{1}\geq l_{1}+l_{2}$.
Hence $l_{3}\geq l_{1}\geq l_{2}\geq l_{3}$ which implies $l_{i}=s/3$
and $2s/3\leq k_{1}$. Furthermore, $s/3\geq k_{2}\geq s/3$
and hence $k_{2}=s/3$ and $k_{1}=2s/3$. \qed

The next Theorem is one of the main results of this section. We will use all the information we obtained about how various conditions on two sequences $\alpha, \beta \in \mathcal{I}_s$ affect the term $T_{\alpha, \beta}$ to obtain a bound on the relation type of $I$ as well as a description of the defining equations of the Rees algebra. 

\begin{theorem}\label {u gens}
Let $R$ be a polynomial ring over a field and let $I$ be a square-free monomial ideal generated by $n$ square-free monomials in $R$. When $n \leq5$, then $\mathcal{R}(I)=S/J$, where $S=R[T_1, \ldots, T_n]$ and $J=SJ_1 + S \cdot (\bigcup \limits_{i=2}^{n-2}J_i)$. In particular, ${\rm rt}(I) \leq n-2$. 
\end{theorem}

\proof
By Theorem~\ref{Taylor}, it suffices  to show that given two sequences $\alpha, \beta \in \mathcal{I}_s$ of length $s>n-2$ then $T_{\alpha, \beta} \in SJ_1 + S \cdot (\bigcup \limits_{i=2}^{n-2}J_i)$. Suppose that there are $l_i$ distinct copies of $a_i$ in $\alpha$ and $k_i$ distinct copies of $b_i$ in $\beta$ and write  $\alpha=(a_{1}^{l_{1}},\ldots,a_{m}^{l_{m}})$,
and $\beta=(b_{1}^{k_{1}},\ldots,b_{t}^{k_{t}})$.

By Lemma~\ref{L3} we may assume $a_{i}\neq b_{j}$
for all $i$,$j$. Also by Lemma~\ref{L5}, we may
assume $1<m \leq 3$ and $1<t \leq 2$, since $n\leq 5$. We are now left with the following cases:
\begin{enumerate}[$($i$)$]
\item Suppose that $\alpha=(a_{1}^{l_{1}},a_{2}^{l_{2}})$ and $\beta=(b_{1}^{k_{1}},b_{2}^{k_{2}})$, where $l_{1}+l_{2}=k_{1}+k_{2}\geq3$. The result follows from Lemma~\ref{L6}.

\item Suppose that $\alpha=(a_{1}^{l_{1}},a_{2}^{l_{2}},a_{3}^{l_{3}})$ and $\beta=(b_{1}^{k_{1}},b_{2}^{k_{2}})$, where $l_{1}+l_{2}+l_{3}=k_{1}+k_{2}\geq4$. The result follows from
Lemma \ref{5gensub}. \qed

\end{enumerate}

\bigskip

The following  example provides a class of square-free monomial ideals  generated by $n>4$
square-free monomials for which the relation type is at least $2n-7$. Notice when $n=5$ one has $2n-7=n-2$. In particular, this establishes that the bound given in Theorem~\ref{u gens} is sharp. We also note that the ideals in the following example are not of fiber type.

\begin{example}\label{rel type ex}{\rm
Let $R=k[y,x_{2},\ldots,x_{n-2},z,u_{2},\ldots,  u_{n-2}]$ be a polynomial ring
over a field $k$ with $n$ a positive integer such that $n>4$. Let $I=(f_{1},\ldots,f_{n})$,
where $f_{1}=(\prod \limits_{i=2}^{n-2}x_{i})z$, $f_{i}=x_{i}y(\prod \limits_{j=2,j\neq i}^{n-2}u_{j})$
for all  $i=2,\ldots, n-2$, $f_{n-1}=(\prod \limits_{i=2}^{n-2}x_{i})y$, and $f_{n}=(\prod \limits_{i=2}^{n-2}u_{i})z$.
Consider \[
F=T_{1}^{n-4}\prod \limits_{i=2}^{n-2}T_{i}-T_{n-1}^{n-3}T_{n}^{n-4} \mbox{ and }
G=zT_{1}^{n-5}\prod \limits_{i=2}^{n-2}T_{i}-yT_{n-1}^{n-5}T_{n}.\]
 It is clear that $F$ and $G$ are in the defining ideal of the Rees algebra of $I$. Moreover, $F$ and $G$ are irreducible and thus $F\in SJ_{2n-7}\setminus S\cdot(\bigcup\limits_{i=1}^{2n-8}J_{i})$
and $G\in SJ_{2n-8}\setminus S\cdot (\bigcup\limits_{i=1}^{2n-9}J_{i})$. Hence the relation type of $I$ is at least $2n-7$. Finally, we note that  $G$ is not in the defining ideal of the special fiber of $I$ and therefore $I$ is not an ideal of fiber type. }
\end{example}

The next lemma establishes conditions for when various generators of the defining ideal of the Rees algebra as in Theorem~\ref{u gens} are irredundant. 

\begin{lemma}\label{5 gens redundant}
Let $R$ be a polynomial ring over a field and let $I$ be a square-free monomial ideal in $R$. 
Let Let $\alpha$, $\beta \in \mathcal{I}_s$  be two sequences of length $s \geq 4$, where $\mathcal{I}_s$ is as in Definition~\ref{sequences}. Suppose that  $\alpha=(a_{1},\ldots,a_{s})$, $\beta=(b_{1},\ldots,b_{1},b_{2})$, where $a_{i} \neq b_j$ for all $i, j$. Suppose that  for all  $i$ there exists a variable  $x_{i} \in R$ such
that $x_{i}\mid f_{a_{j}}$ for all $j\neq i$ and $x_{i}\mid f_{b_{1}}$, $x_{i}\nmid f_{a_{i}}$, and $x_{i}\nmid f_{b_{2}}$. Furthermore, suppose that 
for all $i$ there exists a variable $z_{i}\in R$ such that $z_{i}\mid\gcd(f_{a_{i}},f_{b_{2}})$, $z_{i}\nmid f_{a_{j}}$ for all $j\neq i$ and $z_{i}\nmid f_{b_{1}}$.
Then $T_{\alpha,\beta}\in SJ_{s}\setminus S\cdot(\bigcup \limits_{i=1}^{s-1}J_{i})$.

\end{lemma} 
\proof  We write $f_{a_{i}}=z_{i}h_{i}(\prod \limits_{j=1,j\neq i}^{s}x_{j})$
for all $i=1,\ldots,s$, $f_{b_{1}}=k_{1}(\prod \limits_{i=1}^{s}x_{i})$, and $f_{b_{2}}=k_{2}(\prod \limits_{i=1}^{s}z_{i})$,
for some square-free monomials $h_1,\ldots, h_s, k_1, k_2 \in R$ such that $x_{i}\nmid h_{j}$, $x_{i}\nmid k_{l}$, $z_{i}\nmid h_{j}$,
and $z_{i}\nmid k_{l}$ for all $i$, $j$, $l$. Then we have $f_{\alpha}=\prod \limits_{i=1}^{s}x_{i}^{s-1}\prod \limits_{i=1}^{s}z_{i}\prod \limits_{i=1}^{s}h_{i}$
and $f_{\beta}=k_{1}^{s-1}k_{2}(\prod \limits_{i=1}^{s}x_{i}^{s-1}\prod \limits_{i=1}^{s}z_{i})$.
Let \[
M=\ds{\frac{f_{\alpha}}{\gcd(f_{\alpha},f_{\beta})}=\frac{\prod \limits_{i=1}^{s}h_{i}}{\gcd(\prod \limits_{i=1}^{s}h_{i},k_{1}^{s-1}k_{2})}} \mbox{ and }  
N=\ds{\frac{f_{\beta}}{\gcd(f_{\alpha},f_{\beta})}=\frac{k_{1}^{s-1}k_{2}}{\gcd(\prod \limits_{i=1}^{s}h_{i},k_{1}^{s-1}k_{2})}}.\]
One can observe immediately that $x_{i}\nmid M$ and $z_{i}\nmid N$ for all $i$. Let
$\alpha'$ and $\beta'$ be any proper subsequences of $\alpha$ and $\beta$. Then there exists either  $x_{i}$ or  $z_{j}$ such that $x_{i}\mid \ds{\frac{f_{\alpha'}}{\gcd(f_{\alpha'},f_{\beta'})}}$
or $z_{j}\mid \ds{\frac{f_{\alpha'}}{\gcd(f_{\alpha'},f_{\beta'})}}$ by our assumptions.
Therefore $\ds{\frac{f_{\alpha'}}{\gcd(f_{\alpha'},f_{\beta'})}}\nmid M$, and
similarly, $\ds{\frac{f_{\beta'}}{\gcd(f_{\alpha'},f_{\beta'})}}\nmid N$. This shows that $T_{\alpha,\beta}\in SJ_{s}\setminus S\cdot(\bigcup \limits_{i=1}^{s-1}J_{i})$. \qed

The following Theorem that establishes precisely the defining equations for the Rees algebra of a square-free monomial ideal with up to five generators.

\begin{theorem} \label{5 gens linear}
Let $R$ be a polynomial ring over a field and let $I$ be a square-free monomial ideal generated by $n$ square-free monomials in $R$. Suppose that 
$n \leq 5$. Suppose that there does not exist a pair of sequences $\alpha,\beta \in \mathcal{I}_s$  of length
$s\geq2$ as in Definition~\ref{sequences}, such that the conditions of Lemma~\ref{5 gens redundant} are satisfied. Then  $I$ is an ideal of linear type.
\end{theorem}

\proof Consider a pair of sequences  $\alpha,\beta$
 as in Definition~\ref{sequences} of length $s\geq 2$ such that $\alpha=(a_{1},\ldots,a_{s})$, $\beta=(b_{1},\ldots,b_{1},b_{2})$, where $a_{i} \neq b_j$ for all $i, j$. Suppose that one of the following conditions is not satisfied:
for all  $i$, there exists a variable $x_{i} \in R$ such
that $x_{i}\mid f_{a_{j}}$ for all $j\neq i$ and $x_{i}\mid f_{b_{1}}$, $x_{i}\nmid f_{a_{i}}$, and $x_{i}\nmid f_{b_{2}}$, and for all $i$ there exists a variable $z_{i} \in R$ such that $z_{i}\mid\gcd(f_{a_{i}},f_{b_{2}})$, $z_{i}\nmid f_{a_{j}}$ for all $j\neq i$ and $z_{i}\nmid f_{b_{1}}$.
We will show that  $T_{\alpha,\beta}\in SJ_{1}+S\cdot(\bigcup \limits_{i=1}^{s-1}J_{i})$ and by induction $T_{\alpha, \beta} \in SJ_1$.

By the Theorem~\ref{u gens} and Lemmas~\ref{L3},~\ref{L5} it is enough to consider sequences
of length 2 or 3, i.e. either $\alpha=(a_{1},a_{2})$ and $\beta=(b_{1},b_{2})$ or $\alpha=(a_{1},a_{2},a_{3})$ and $\beta=(b_{1},b_{1},b_{2})$
. 
The proof for
the first case can be treated as a special case of the second and thus we will only consider the case with $\alpha=(a_{1},a_{2},a_{3})$ and $\beta=(b_{1},b_{1},b_{2})$. Without lost of generality, we will
show that if there does not exist a variable $x\in R$ such that $x\mid\gcd(f_{a_{2}},f_{a_{3}})$, $x\mid f_{b_{1}}$, $x\nmid f_{a_{1}}$, and $x\nmid f_{b_{2}}$
then \begin{eqnarray*}
\gcd(f_{a_{1}}f_{a_{2}}f_{a_{3}},f_{b_{1}}^{2}f_{b_{2}}) & = & \frac{\gcd(f_{a_{1}},f_{b_{1}})\gcd(f_{\alpha},f_{b_{1}}f_{b_{2}})}{A}\\
 & = & \frac{\gcd(f_{a_{1}},f_{\beta})\gcd(f_{a_{2}}f_{a_{3}},f_{b_{1}}f_{b_{2}})}{B},\end{eqnarray*}
for some $A,B \in R$. 

Similarly, if there does not exist a variable $z\in R$ such that $z\mid\gcd(f_{a_{1}},f_{b_{2}})$, $z\nmid f_{a_{j}}$ for $j=2,3$, and $z\nmid f_{b_{1}}$ then 

\begin{eqnarray*}
\gcd(f_{a_{1}}f_{a_{2}}f_{a_{3}},f_{b_{1}}^{2}f_{b_{2}}) & = & \frac{\gcd(f_{a_{1}},f_{b_{1}})\gcd(f_{a_{2}}f_{a_{3}},f_{\beta})}{C}\\
 & = & \frac{\gcd(f_{\alpha},f_{b_{1}})\gcd(f_{a_{2}}f_{a_{3}},f_{b_{1}}f_{b_{2}})}{D},\end{eqnarray*}
 for some $C,D \in R$. Hence we may apply Lemma~\ref{Tab split extended} to obtain the result.

 For the first part, we notice that
\begin{eqnarray*}
\gcd(f_{a_{1}}f_{a_{2}}f_{a_{3}},f_{b_{1}}^{2}f_{b_{2}}) & = & \frac{\gcd(f_{\alpha},f_{b_{1}})\gcd(f_{\alpha},f_{b_{1}}f_{b_{2}})}{E}\\
 & = & \frac{\gcd(f_{a_{1}},f_{b_{1}})\gcd(f_{a_{2}}f_{a_{3}},f_{b_{1}})\gcd(f_{\alpha},f_{b_{1}}f_{b_{2}})}{EF},\end{eqnarray*}
\begin{eqnarray*}
\gcd(f_{a_{1}}f_{a_{2}}f_{a_{3}},f_{b_{1}}^{2}f_{b_{2}}) & = & \frac{\gcd(f_{a_{1}},f_{\beta})\gcd(f_{a_{2}}f_{a_{3}},f_{\beta})}{G}\\
 & = & \frac{\gcd(f_{a_{1}},f_{\beta})\gcd(f_{a_{2}}f_{a_{3}},f_{b_{1}})\gcd(f_{a_{2}}f_{a_{3}},f_{b_{1}}f_{b_{2}})}{GH},\end{eqnarray*}
for some $E, F, G, H \in R$. It is enough to show
for any variable $x \in R$ if $x\mid\gcd(f_{a_{2}}f_{a_{3}},f_{b_{1}})$, then $x\mid EF$
and $x\mid GH$. This is equivalent to showing that  for any integer $u$
if $x^{u}\mid\gcd(f_{\alpha},f_{\beta})$ then $x^{u}|\gcd(f_{a_{1}},f_{b_{1}})\gcd(f_{\alpha},f_{b_{1}}f_{b_{2}})$
and $x^{u}\mid\gcd(f_{a_{1}},f_{\beta})\gcd(f_{a_{2}}f_{a_{3}},f_{b_{1}}f_{b_{2}})$.
This is straightforward to verify and the only case that is not trivial
is when  $x\nmid f_{a_{1}}$ and $x\mid\gcd(f_{a_{2}},f_{a_{3}})$. Then $x\mid f_{b_{2}}$ by assumption. Hence $x^{2}\mid\gcd(f_{\alpha},f_{b_{1}}f_{b_{2}})$
and $x^{2}\mid\gcd(f_{a_{2}}f_{a_{3}},f_{b_{1}}f_{b_{2}})$.

We can use a similar argument for the second part and the only case 
that is not trivial is when $x\mid\gcd(f_{a_{1}},f_{b_{1}}f_{b_{2}})$
and $x\nmid f_{b_{1}}$. Then $x\mid f_{b_{2}}$ and $x\mid f_{a_{2}}$
or $x\mid f_{a_{3}}$ by assumption. Hence $x^{2}\nmid\gcd(f_{\alpha},f_{\beta})$,
$x\mid\gcd(f_{a_{2}}f_{a_{3}},f_{\beta})$, and $x\mid\gcd(f_{a_{2}}f_{a_{3}},f_{b_{1}}f_{b_{2}})$. \qed

We now introduce a construction that is also known as the $1$-skeleton or the line graph of the hypergraph of a monomial ideal.

\begin{construction}\label{graph til}{\rm
Let $R$ be a polynomial ring over a field and let $I$ be a monomial ideal  in $R$ . Let $f_1, \ldots, f_n$ be a minimal monomial generating set of $I$. 
We construct the following graph $\til{G}(I)$:
For each $f_i$  we associate a vertex $y_i$ to it. The edges of this graph are $\{y_i, y_j\}$, where ${\gcd}(f_i, f_j)\neq 1$. We call the graph $\til{G}(I)$ the \it{generator graph} of $I$.}
\end{construction}

The purpose of introducing the generator graph of a monomial ideal is to utilize the graph structure in order to determine combinatorial conditions that determine the defining equations of the Rees algebra.  We first observe the following.

\begin{remark}\label{even closed walk}{\rm 
With the same assumptions as in Lemma~\ref{5 gens redundant} the sequences $\alpha, \beta$ correspond to an even closed walk in $\til{G}(I)$. Indeed, the following \begin{eqnarray*}
 & \{y_{b_{1}},\{y_{b_{1}},y_{a_{1}}\},y_{a_{1}},\{y_{a_{1}},y_{b_{1}}\},y_{b_{1}},\{y_{b_{1}},y_{a_{2}}\},y_{a_{2}},\{y_{a_{2}},y_{b_{1}}\},y_{b_{1}},\ldots,\\
   & y_{b_{1}},\{y_{b_{1}},y_{a_{s-1}}\},y_{a_{s-1}},\{y_{a_{s-1}},y_{b_{2}}\},y_{b_{2}},\{y_{b_{2}},y_{a_{s}}\},y_{a_{s}},\{y_{a_{s}},y_{b_{1}}\},y_{b_{1}}\}\end{eqnarray*}
 is an even closed walk and $T_{\alpha,\beta}\in SJ_{s}\setminus S\cdot(\bigcup \limits_{i=1}^{s-1}J_{i})$.
}

\end{remark}

We now introduce the notion of a subgraph induced by two sequences.

 \begin{definition} \label{induced graph}{\rm
 Let $R$ be a polynomial ring over a field and let $I$ be a monomial ideal  in $R$.
 Let $\alpha$, $\beta \in \mathcal{I}_s$  be two sequences of length $s \geq 2$, where $\mathcal{I}_s$ is as in Definition~\ref{sequences}.  
 Let $\alpha=(a_{1}^{l_{1}},a_{2}^{l_{2}},\ldots,a_{m}^{l_{m}})$
and $\beta=(b_{1}^{r_{1}},\ldots,b_{t}^{r_{t}})$, where $l_i$ and $r_i$ are the number of distinct copies of $a_i$ and $b_i$, respectively. Let $\alpha_1=(a_1, \ldots, a_m)$ and $\beta_1=(b_{i_1},\ldots, b_{i_r})$, where $\{b_{i_1}, \ldots, b_{i_r}\}=\{b_1, \ldots, b_t\} \setminus\{a_1, \ldots, a_m\}$. The last condition ensures that $a_i \neq b_{i_k}$ for all $i,i_{k}$. 
Consider $K=(f_{a_{1}}, \ldots, f_{a_{m}}, f_{b_{i_1}}, \ldots, f_{b_{i_r}})$ and the generator graph $\til{G}(K)$ of $K$. Notice that $K \subset I$ and that $\til{G}(K) $ is a subgraph of $\til{G}(I)$, since $a_i \neq b_{i_k}$ for all $i,i_{k}$.  We call $\til{G}(K)$ the graph {\it induced by} $\alpha$ and $\beta$.}
 \end{definition}

\begin{remark} \label{disconnected}{\rm
Let $R$ be a polynomial ring over a field and let $I$ be a monomial ideal in $R$. Let $\alpha, \beta \in \mathcal{I}_s$ be two sequences of length $s \geq 2$,  where $\mathcal{I}_s$ is as in Definition~\ref{sequences}. 
Let $G$ be the graph induced by $\alpha$ and $\beta$.  If $G$ is a disconnected graph then $T_{\alpha, \beta} \in SJ_1+S\cdot (\bigcup \limits_{i=2}^{s-1}J_i)$, by  Proposition~\ref{L2}. }
\end{remark}

The following Theorem is an extension of Theorem~\ref{u gens}. We use the generator graph to give a description for the defining equations of the Rees algebra.

\begin{theorem} \label{what survives}
Let $R$ be a polynomial ring over a field and let $I$ be a square-free monomial ideal generated by $n$ square-free monomials in $R$. Suppose that the generator graph of $I$ is the disjoint union of graphs with at most $5$ vertices. Let $\mathcal{R}(I)=S/J$, where $S=R[T_1, \ldots, T_n]$. Then $J=SJ_1 + S \cdot (J_2\cup J_3)$ and in particular, ${\rm rt}(I) \leq 3$. Furthermore, if  $T_{\alpha, \beta} \in S \cdot (J_2\cup J_3) \setminus SJ_1$, then the subgraph induced by $\alpha, \beta$ is an even closed walk of $\til{G}(I)$.

\end{theorem}

\proof  By Remark~\ref{disconnected} we may assume that $\til{G}(I)$ is connected. Then
the first part follows immediately by Theorem~\ref{u gens}. The last part follows from Theorem~\ref{5 gens linear}, Lemma~\ref{5 gens redundant}, and Remark~\ref{even closed walk}. \qed

The next example illustrates the result of Theorem~\ref{what survives}.

\begin{example}{\rm
Let $R=k[x_1, \ldots, x_7]$ be a polynomial ring over a field $k$. Let $I=(f_1, \ldots, f_5)$, where $f_1=x_{1}x_{2}x_{3}$, $f_2=x_{1}x_{2}x_{4}x_{7}$, $f_3=x_{2}x_{3}x_{6}$, $f_4=x_{4}x_{5}x_{6}$, and $f_5=x_{1}x_{3}x_{5}$. Let $\rees{I}=S/J$, where $S=R[T_1,\ldots, T_5]$. 
Let $\alpha=(1,1,4)$ and $\beta=(2,3,5)$. The subgraph induced by $\alpha, \beta$ is the even closed walk $$\{y_1, \{y_1,y_2\},y_2, \{y_1,y_2\},y_1, \{y_1,y_3\},y_3, \{y_3, y_4\},y_4, \{y_4, y_5\},y_5,\{y_1,y_5\}\}$$ of the generator graph of $I$ shown below:

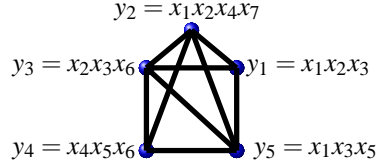
\begin{center}
\begin{tikzpicture}
\shade [shading=ball, ball color=blue] (0,1) circle (0.1) node [above ] {\textsf{\textbf{\textcolor{black}{\small{$y_2=x_{1}x_{2}x_{4}x_{7}$ }}}}};
\shade [shading=ball, ball color=blue] (0.6,0.5) circle (0.1) node [right] {\textsf{\textbf{\textcolor{black}{\small{$y_1=x_{1}x_{2}x_{3}$}}}}}; 
\shade [shading=ball, ball color=blue] (0.6,-0.6) circle (0.1) node [right ]  {\textsf{\textbf{\textcolor{black}{\small{ $y_5=x_{1}x_{3}x_{5}$}}}}};
\shade [shading=ball, ball color=blue] (-0.6,-0.6) circle (0.1) node [left ]  {\textsf{\textbf{\textcolor{black}{\small{ $y_4=x_{4}x_{5}x_{6}$}}}}};
\shade [shading=ball, ball color=blue] (-0.6,0.5) circle (0.1) node [left ]  {\textsf{\textbf{\textcolor{black}{ \small{$y_3=x_{2}x_{3}x_{6}$}}}}};

\draw [line width=2pt  ] (0,1)--(0.6,0.5);
\draw [line width=2pt  ] (0.6,0.5)--(0.6,-0.6);
\draw [line width=2pt  ] (0.6,-0.6)--(-0.6,-0.6);
\draw [line width=2pt  ] (-0.6,-0.6)--(-0.6,0.5);
\draw [line width=2pt  ] (0,1)--(-0.6,-0.6);
\draw [line width=2pt  ] (-0.6,0.5)--(0.6,0.5);
\draw [line width=2pt  ] (0,1)--(0.6,-0.6);
\draw [line width=2pt  ] (0.6,-0.6)--(-0.6,0.5);
\draw [line width=2pt  ] (0,1)--(-0.6,0.5);

\end{tikzpicture}

\end{center}
Then by Theorem~\ref{what survives} we have that $T_{\alpha, \beta}= T_2T_3T_5- x_7T_1^2T_4  \in SJ_3 \setminus (SJ_1+SJ_2)$, since it is irreducible.
}
\end{example}

The following example comes from  \cite[Example~3.1]{Vill}. Villarreal uses this example to show that his methods do not extend in the case of square-free monomial ideals generated in  degree higher than 2. In light of Theorem~\ref{what survives} we now give an explicit description of the defining equations of the Rees algebra for this example.

\begin{example}\label{Vil ex} {\rm
Let $R=k[x_{1},\ldots, x_{7}]$, where $k$ is a field and let $I$ be an ideal of $R$ generated
by $f_{1}=x_{1}x_{2}x_{3},$ $f_{2}=x_{2}x_{4}x_{5}$, $f_{3}=x_{5}x_{6}x_{7}$,
$f_{4}=x_{3}x_{6}x_{7}$. 
Then the defining ideal of $\mathcal{R}(I)$
is minimally generated by the binomials:\[
x_{3}T_{3}-x_{5}T_{4},x_{6}x_{7}T_{1}-x_{1}x_{2}T_{4},x_{6}x_{7}T_{2}-x_{2}x_{4}T_{3},x_{4}x_{5}T_{1}-x_{1}x_{3}T_{2},x_{4}T_{1}T_{3}-x_{1}T_{2}T_{4}.\]
Notice that the only binomial that is not linear is $x_{4}T_{1}T_{3}-x_{1}T_{2}T_{4}$
and it comes from the unique even cycle of $\widetilde{G}(I)$, which in this case
is the graph of a square. }
\end{example}

\section{Square-free monomial ideals of linear type}\label{lin type}
 
 In this section we turn our attention to ideals of linear type. Villarreal showed that when $I$ is the edge ideal of a  graph then $I$ is of linear type if and only if the graph of $I$  is the disjoint union of graphs of trees or graphs that have a unique odd cycle, \cite[Corollary~3.2]{Vill}. Inspired by this result we use the generator graph associated to any monomial ideal as in Construction~\ref{graph til}  in order to obtain similar results as Villarreal for any square-free monomial ideal without restrictions on the degrees of the generators. A natural first class to consider is ideals with generator graph the graph of a forest.

\begin{proposition}\label{tree linear type}
 Let $R$ be a polynomial ring over a field and let $I$ be a square-free monomial ideal in $R$. We further assume that the generator graph $\til{G}(I)$ is the graph of a forest. Then $I$ is an ideal of linear type.
\end{proposition}

\proof Let $f_1, \ldots, f_n$ be a minimal monomial generating set for $I$.
We show that  for all sequences $\alpha, \beta \in \mathcal{I}_s$ of length $s\geq 2$, where $\mathcal{I}_s$ is as in Definition~\ref{sequences}, we have $T_{\alpha, \beta} \in SJ_1+S\cdot (\bigcup \limits_{i=2}^{s-1}J_i)$ and by induction it follows that $T_{\alpha, \beta} \in SJ_1$.

Since $\til{G}(I)$ is the graph of a forest then every subgraph of $\til{G}(I)$ is also a forest. By Proposition~\ref{L2} we may assume that $\til{G}(I)$ is connected, i.e. it is the graph of a tree. Notice that for every graph of a tree there exists a vertex that is only connected to one other vertex. 

Suppose that $s \geq 2$. Let $\alpha=(a_{1}^{l_{1}},a_{2}^{l_{2}},\ldots,a_{m}^{l_{m}})$
and $\beta=(b_{1}^{r_{1}},\ldots,b_{t}^{r_{t}})$, where $l_i$ and $r_i$ are the number of distinct copies of $a_i$ and $b_i$, respectively. 
Notice that if $a_i =b_j$ for some $i$ and some $j$ then the result follows by Lemma~\ref{L3} and the induction hypothesis. Hence we may assume that $a_i \neq b_j$ for all $i,j$.
Let $\alpha_1=(a_1, \ldots, a_m)$ and $\beta_1=(b_1,\ldots, b_t)$. Then the graph $G'$ induced by $\alpha_1$, $\beta_1$  is a subgraph of $\til{G}(I)$, by Definition~\ref{induced graph}. Hence $G'$ is  the graph of a forest. If $G'$ is disconnected then the result follows by Proposition~\ref{L2} and the induction hypothesis. Therefore, we may assume that $G'$ is connected and hence the graph of  a tree. Thus without loss of generality we may assume that $y_{a_1}$ is only connected to $y_{b_1}$.
 
If $l_{1}\geq r_{1}$, then $\gcd(f_{a_{1}},f_{b_{2}}^{r_{2}}\cdots f_{b_{t}}^{r_{t}})=1$ and hence  $T_{\alpha, \beta} \in SJ_1+S\cdot (\bigcup \limits_{i=2}^{s-1}J_i)$, by Proposition~\ref{L2}. So we may assume that  $l_{1}<r_{1}$. Then we claim  that 
\begin{enumerate}[$($a$)$]
\item $\gcd(f_{\alpha},f_{\beta})=\gcd(f_{a_{1}}^{l_{1}},f_{b_{1}}^{l_{1}})\gcd(f_{a_{2}}^{l_{2}} \cdots f_{a_{m}}^{l_{m}},f_{\beta})$,
\item $\gcd(f_{\alpha},f_{\beta})=\ds{\frac{\gcd(f_{\alpha},f_{b_{1}}^{l_{1}})\gcd(f_{a_{2}}^{l_{2}}\cdots f_{a_{m}}^{l_{k}},f_{b_{1}}^{r_{1}-l_{1}}f_{b_{2}}^{r_{2}}\cdots f_{b_{t}}^{r_{t}})}{A}}$,
\end{enumerate}
for some $A\in R$. Notice that by the claim and Lemma~\ref{Tab split extended} it follows that  $T_{\alpha, \beta} \in SJ_1+S\cdot (\bigcup \limits_{i=2}^{s-1}J_i)$. Therefore, it remains to prove the claim.

For part~(a) notice that
$\gcd(f_{\alpha},f_{\beta})=\gcd(f_{a_{1}}^{l_{1}},f_{\beta})\gcd(f_{a_{2}}^{l_{2}}\cdots f_{a_{m}}^{l_{m}},f_{\beta})=(f_{a_{1}}^{l_{1}},f_{b_{1}}^{l_{1}})\gcd(f_{a_{2}}^{l_{2}}\cdots f_{a_{m}}^{l_{m}},f_{\beta})$,
where the equalities follow from the fact that $\gcd(f_{a_{1}},f_{a_{i}})=\gcd(f_{a_{1}},f_{b_{j}})=1$
for any $i>1, j>1$ and Remark~\ref{gcd properties}~(b).

For part (b), notice that \begin{eqnarray*}
\gcd(f_{\alpha},f_{\beta}) & = & \gcd(f_{\alpha},f_{b_{1}}^{l_{1}})\gcd(f_{\alpha},f_{b_{1}}^{r_{1}-l_{1}}f_{b_{2}}^{r_{2}}\cdots f_{b_{t}}^{r_{t}})/B\\
 & = & \gcd(f_{\alpha},f_{b_{1}}^{l_{1}})\gcd(f_{a_{1}}^{l_{1}},f_{b_{1}}^{r_{1}-l_{1}}f_{b_{2}}^{r_{2}}\cdots f_{b_{t}}^{r_{t}})\gcd(f_{a_{2}}^{l_{2}}\cdots f_{a_{m}}^{l_{m}},f_{b_{1}}^{r_{1}-l_{1}}f_{b_{2}}^{r_{2}}\cdots f_{b_{t}}^{r_{t}})/BC\\
 & = & \gcd(f_{\alpha},f_{b_{1}}^{l_{1}})\gcd(f_{a_{1}}^{l_{1}},f_{b_{1}}^{r_{1}-l_{1}})\gcd(f_{a_{2}}^{l_{2}}\cdots f_{a_{m}}^{l_{m}},f_{b_{2}}^{r_{2}}\cdots f_{b_{t}}^{r_{t}})/BC,\end{eqnarray*}
for some $B,C \in R$, since $\gcd(f_{a_{1}}^{l_{1}}, f_{b_{2}}^{r_{2}}\cdots f_{b_{t}}^{r_{t}})=1$. It is enough to show that $\gcd(f_{a_{1}}^{l_{1}},f_{b_{1}}^{r_{1}-l_{1}})\mid BC$,
i.e. for any variable $x \in R$ and any integer $v$ if $ $$x^{v}\mid\gcd(f_{a_{1}}^{l_{1}},f_{b_{1}}^{r_{1}-l_{1}})$,
then $x^{v}\mid BC$. This is equivalent to showing that  for any variable $x \in R$ if $x\mid\gcd(f_{a_{1}}^{l_{1}},f_{b_{1}}^{r_{1}-l_{1}})$
and $x^{u}\mid\gcd(f_{\alpha},f_{\beta})$ for some integer $u$, then
$x^{u}\mid\gcd(f_{\alpha},f_{b_{1}}^{l_{1}})\gcd(f_{a_{2}}^{l_{2}}\cdots f_{a_{m}}^{l_{m}},f_{b_{1}}^{r_{1}-l_{1}}f_{b_{2}}^{r_{2}}\cdots f_{b_{t}}^{r_{t}})$.
If $x\mid\gcd(f_{a_{1}}^{l_{1}},f_{b_{1}}^{r_{1}-l_{1}})$ and $x^{u}\mid\gcd(f_{\alpha},f_{\beta})$,
then $u\leq l_{1}$ by the fact that $y_{a_{1}}$ is only connected
to $y_{b_{1}}$ and $l_{1}<r_{1}$. Therefore $x^{u}\mid\gcd(f_{\alpha},f_{b_{1}}^{l_{1}})$.\qed

\bigskip

The following lemma allows us to handle the induction step in the case of odd cycles, in order to prove that when the generator graph of a square-free monomial ideal is the graph of an odd cycle, then the ideal is of linear type.

\begin{lemma} \label{odd cycle lemma}
 Let $R$ be a polynomial ring over a field and let $I$ be a square-free monomial ideal in $R$.  Let $\alpha$, $\beta \in \mathcal{I}_s$  be two sequences of length $s \geq 4$, where $\mathcal{I}_s$ is as in Definition~\ref{sequences}. 
 Let $\alpha=(a_{1}^{l_{1}},a_{2}^{l_{2}},\ldots,a_{m}^{l_{m}})$
and $\beta=(b_{1}^{r_{1}},\ldots,b_{t}^{r_{t}})$,  where $l_i$ and $r_i$ are the number of distinct copies of $a_i$ and $b_i$, respectively. Suppose that $l_{1}<r_{1}$, $l_{2}<r_{2}$, and  that the graph $\til{G}(I)$ is the graph of an odd cycle of length at least $5$. We further assume that $y_{a_{1}}$ is connected to $y_{b_{1}}$ and $y_{a_{2}}$ only, and  $y_{a_{2}}$ is connected
to $y_{a_{1}}$ and $y_{b_{2}}$ only.  Then $T_{\alpha, \beta} \in SJ_1 + S \cdot (\bigcup \limits_{i=2}^{s-1}J_i)$.

\end{lemma}

\proof 
Let $\alpha'=(a_{3}^{l_{3}},\ldots,a_{m}^{l_{m}})$ and $\beta'=(b_{1}^{r_{1}-l_{1}},b_{2}^{r_{2}-l_{2}},b_{3}^{r_{3}}, \ldots, b_{t}^{r_t})$. Using Remark~\ref{gcd properties} and Remark~\ref{more gcd prop} we have that 
\begin{eqnarray*}\gcd(f_{\alpha},f_{\beta})&=&\gcd(f_{a_{1}}^{l_{1}}f_{a_{2}}^{l_{2}},f_{\beta})\gcd(f_{\alpha'},f_{\beta})\\
&=&\gcd(f_{a_{1}}^{l_{1}}f_{a_{2}}^{l_{2}},f_{b_{1}}^{l_{1}}f_{b_{2}}^{l_{2}})
\gcd(f_{\alpha'},f_{\beta}),
\end{eqnarray*}
 since $l_{1}<r_{1}$ and
$l_{2}<r_{2}$.

We claim that $\gcd(f_{\alpha}, f_{\beta})=\ds{\frac{\gcd(f_{\alpha},f_{b_{1}}^{l_{1}}f_{b_{2}}^{l_{2}})(f_{\alpha'},f_{\beta'})}{C}}$, for some $C \in R$. Then the result will follow by  Lemma~\ref{Tab split extended}. Notice that \begin{eqnarray*}
\gcd(f_{\alpha},f_{\beta}) & = & \frac{\gcd(f_{\alpha},f_{b_{1}}^{l_{1}}f_{b_{2}}^{l_{2}})\gcd(f_{\alpha},f_{\beta'})}{D}\\
 & = & \frac{\gcd(f_{\alpha},f_{b_{1}}^{l_{1}}f_{b_{2}}^{l_{2}})\gcd(f_{\alpha'},f_{\beta'})\gcd(f_{a_{1}}^{l_{1}}f_{a_{2}}^{l_{2}},f_{\beta'})}{DE}\\
 & = & \frac{\gcd(f_{\alpha},f_{b_{1}}^{l_{1}}f_{b_{2}}^{l_{2}})\gcd(f_{\alpha'},f_{\beta'})\gcd(f_{a_{1}}^{l_{1}}f_{a_{2}}^{l_{2}},f_{b_{1}}^{r_{1}-l_{1}}f_{b_{2}}^{r_{2}-l_{2}})}{DE},\end{eqnarray*}
where $D,E\in R$ and the third equality follows from the fact that
$y_{a_{1}}$ and $y_{a_{2}}$ are only connected to $y_{b_{1}}$ and $y_{b_{2}}$, respectively.
It is enough to show that for any variable $x \in R$ and any integer $v$ if  $x^{v}\mid\gcd(f_{a_{1}}^{l_{1}}f_{a_{2}}^{l_{2}},f_{b_{1}}^{r_{1}-l_{1}}f_{b_{2}}^{r_{2}-l_{2}})$,
then $x^{v}\mid DE$. This is equivalent to showing that  if $x\mid\gcd(f_{a_{1}}^{l_{1}}f_{a_{2}}^{l_{2}},f_{b_{1}}^{r_{1}-l_{1}}f_{b_{2}}^{r_{2}-l_{2}})$
and $x^{u}\mid\gcd(f_{\alpha},f_{\beta})$, then $x^{u}\mid\gcd(f_{\alpha},f_{b_{1}}^{l_{1}}f_{b_{2}}^{l_{2}})\gcd(f_{\alpha'},f_{\beta'})$.
If $ $$x\mid\gcd(f_{a_{1}}^{l_{1}}f_{a_{2}}^{l_{2}},f_{b_{1}}^{r_{1}-l_{1}}f_{b_{2}}^{r_{2}-l_{2}})$,
then $x$ can not divide both $f_{a_{1}}$ and $f_{a_{2}}$, since otherwise  $\{y_{a_{1}},y_{a_{2}},y_{b_{1}}\}$ or $\{y_{a_1}, y_{a_2}, y_{b_{2}}\}$ will form a $3$-cycle in $\til{G}(I)$.
Similarly, $x$ can not divide both $f_{b_{1}}$ and $f_{b_{2}}$. Therefore, without loss of generality we may assume $x\mid f_{a_{1}}$ and
$x\mid f_{b_{1}}$. Since $l_{1}<r_{1}$ and $y_{a_{1}}$ is not connected
to $y_{a_{i}}$ for all $i>2$, we obtain that if $x^{u}\mid\gcd(f_{\alpha},f_{\beta})$,
then $u\leq l_{1}$. Therefore $x^{u}\mid\gcd(f_{\alpha},f_{b_{1}}^{l_{1}}f_{b_{2}}^{l_{2}})$.
 \qed

In the next Proposition we are able to handle the case where each of the connected components of the generator graph has a unique odd cycle.

\begin{proposition}\label{odd cycle linear type}
Let $R$ be a polynomial ring over a field and let $I$ be a square-free monomial ideal in $R$.  We further assume that generator graph $\til{G}(I)$ of $I$ is the disjoint union of graphs 
with a unique odd cycle. Then $I$ is an ideal of linear type.
\end{proposition}

\proof 
By Remark~\ref{disconnected} it suffices to consider the connected components of $\til{G}(I)$. Hence we may assume that $\til{G}(I)$ is connected and it has a unique odd cycle.  We show that for all  sequences $\alpha, \beta \in \mathcal{I}_s$  of length $s \geq 2$, where $\mathcal{I}_s$ is  as in Definition~\ref{sequences}, we have $T_{\alpha, \beta} \in SJ_1+S\cdot (\bigcup \limits_{i=2}^{s-1}J_i)$. By induction it will follow  that   $T_{\alpha, \beta} \in SJ_1$.

 Let $\alpha=(a_{1}^{l_{1}},a_{2}^{l_{2}},\ldots,a_{m}^{l_{m}})$
and $\beta=(b_{1}^{r_{1}},\ldots,b_{t}^{r_{t}})$,  where $l_i$ and $r_i$ are the number of distinct copies of $a_i$ and $b_i$, respectively. Using Lemma~\ref{L3} we may assume without loss of generality that $a_i \neq b_j$ for all $i$ and $j$. 
Let $\alpha_1=(a_1, \ldots, a_m)$ and $\beta_1=(b_1, \ldots, b_t)$. Then the graph $G'$
induced by $\alpha_1$ and $\beta_1$ is a subgraph of $\til{G(I)}$, by Definition~\ref{induced graph}. By Proposition~\ref{L2} we may assume that $G'$ is a connected graph. 
Notice that if  there exists a vertex $y_i$ that is only connected to one of the $y_j$ then $T_{\alpha, \beta} \in SJ_1+S\cdot (\bigcup \limits_{i=2}^{s-1}J_i)$ by the proof of Proposition~\ref{tree linear type}. Hence we may assume that $G'$ is the graph of an odd cycle. Therefore every vertex is only connected to two other vertices.

Let $l=m+t$. Since $G'$ is the graph of an odd cycle then $l$ is odd and hence  $m \neq t$.  We proceed by induction on $l$. If $l =3$ then $\til{G}(I)$ is the graph of a triangle and the result follows from Lemma~\ref{L3} and Lemma~\ref{L5}. Suppose that $l \geq 5$. Using Lemma~\ref{L5} we may assume that $m \geq 2$ and $t \geq 2$. Notice that  $y_{a_1}$ is connected to two other vertices. Without loss of generality we may assume that we  have the following possible cases:
\begin{enumerate}[(i)]
\item $y_{a_1}$ is connected only to $y_{a_2}$ and $y_{a_3}$,
\item $y_{a_1}$ is connected only to $y_{a_2}$ and $y_{b_1}$,
\item $y_{a_1}$ is connected only to $y_{b_1}$ and $y_{b_2}$.
\end{enumerate}
We handle each case separately. 
For case (i) we note that  since $y_{a_1}$ is connected only to $y_{a_2}$ and $y_{a_3}$, then using $k=l_1$ in Proposition~\ref{L2} and induction we obtain the result.

For case (ii) suppose that $y_{a_1}$ is connected only to $y_{a_2}$ and $y_{b_1}$. If $l_1 \geq r_1$ then using $k=l_1$ in Proposition~\ref{L2} and induction we obtain the result. If $l_1<r_1$ then we consider $y_{a_2}$. Notice that $y_{a_2}$ can not be connected to $y_{b_1}$ since then $\{y_{a_1},y_{a_2},y_{b_1}\}$ will form a triangle, which is a contradiction since $l\geq 5$. Hence either $y_{a_2}$ is  connected to $y_{a_3}$, or $y_{a_2}$ is connected to $y_{b_i}$ for some $i \geq 2$. If $y_{a_2}$ is connected to $y_{a_3}$ then case (i) yields the result.

Suppose that $y_{a_2}$ is connected to $y_{b_i}$ for some $i \geq 2$. Without loss of generality suppose that $i=2$. If $l_2 \geq r_2$ then we may use $k=l_2$ in Proposition~\ref{L2} and induction. If $l_2 <r_2$ and since $l_1 <r_1$ we can use Lemma~\ref{odd cycle lemma} and induction.

Finally for case (iii) suppose that $y_{a_1}$ is connected only to $y_{b_1}$ and $y_{b_2}$. If every $y_{a_i}$ is connected to $y_{b_{i_1}}$ and $y_{b_{i_2}}$ for $i_1,  i_2 \neq 1,2$ then either $m=t$, which is impossible since $l$ is odd or there exist $i \neq j$ such that $y_{b_i}$ is connected to $y_{b_j}$. Then by switching the role of $\alpha$ and $\beta$ we are in either case (i) or case (ii). \qed

We are now ready to state the main theorem of this section.

\begin{theorem}\label{union graph linear type}
Let $R$ be a polynomial ring over a field and let $I$ be a square-free monomial ideal in $R$.  We further assume that the generator graph, $\til{G}(I)$, is the disjoint union of graphs of trees and graphs that with a unique odd cycle. Then $I$ is an ideal of linear type.

\end{theorem}

\proof
The result follows immediately by Proposition~\ref{tree linear type}, Proposition~\ref{odd cycle linear type}, and Proposition~\ref{L2}. \qed

We conclude this article with the following remark.

\begin{remark}{\rm
Let $R$ be a polynomial ring over a field and let $I$ be a square-free monomial ideal generated by $3$ square-free monomials in $R$. Then the fact that $I$ is of linear type is already known, but it also follows immediately from Proposition~\ref{tree linear type} and Proposition~\ref{odd cycle linear type}.

 }
\end{remark}


\end{document}